\documentclass[12pt,twoside]{article}
\usepackage{amsmath}
\usepackage[all]{xy}

\newcommand{\VolumeHeader}{}
\newcommand{\VolumeSerial}{LNS}

\newcommand{\ActivityName}{ {\normalsize {\it 
Summer School on High-dimensional Manifold Topology }}}
\newcommand{\ActivityDate}{ {\normalsize {\it
Trieste, 21 May -- 8 June 2001
}}}

\usepackage{amssymb}
\newcommand{\be}{\begin{equation}}
\newcommand{\ee}{\end{equation}}
\newcommand{\bea}{\begin{eqnarray}}
\newcommand{\eea}{\end{eqnarray}}
\def\A{\mathbb A}

\def\K{\mathbb K}
\def\L{\mathbb L}
\def\Q{\mathbb Q}

\def\S{\mathbb S}
\def\Z{\mathbb Z}
\def\qed{$\square$}
\newcommand{\LectureHeader}{The structure set}

\begin{document}
\pagestyle{myheadings}
\markboth{\LectureHeader}{\VolumeHeader}
\markright{\VolumeHeader}


\begin{titlepage}


\title{The structure set of an arbitrary space,\\
the algebraic surgery exact sequence\\
and the total surgery obstruction}

\author{Andrew Ranicki\thanks{aar@maths.ed.ac.uk}
\\[1cm]
{\normalsize
{\it Department of Mathematics and Statistics}}\\
{\normalsize
{\it University of Edinburgh, Scotland, UK} }
\\[9cm]
{\normalsize {\it Lecture given at the: }}
\\
\ActivityName 
\\
\ActivityDate 
\\[1cm]
{\small \VolumeSerial} 
}
\date{}
\maketitle
\thispagestyle{empty}
\end{titlepage}

\baselineskip=14pt
\newpage
\thispagestyle{empty}


\begin{abstract}

The algebraic theory of surgery gives a necessary and sufficient chain
level condition for a space with $n$-dimensional Poincar\'e duality to
be homotopy equivalent to an $n$-dimensional topological manifold.
A relative version gives a necessary and sufficient
chain level condition for a simple homotopy equivalence of 
$n$-dimensional topological manifolds to be homotopic to a homeomorphism.
The chain level obstructions come from a chain level interpretation
of the fibre of the assembly map in surgery. 

The assembly map $A:H_n(X;\L_{\bullet}) \to L_n(\Z[\pi_1(X)])$ is a
natural transformation from the generalized homology groups of a space
$X$ with coefficients in the 1-connective simply-connected surgery
spectrum $\L_{\bullet}$ to the non-simply-connected surgery obstruction
groups $L_*(\Z[\pi_1(X)])$. The $(\Z,X)$-category has objects
based f.g. free $\Z$-modules with an $X$-local structure.
The assembly maps $A$ are induced by a functor from the
$(\Z,X)$-category to the category of based f.g. free $\Z[\pi_1(X)]$-modules.
The generalized homology groups $H_*(X;\L_{\bullet})$ are the cobordism
groups of quadratic Poincar\'e complexes over $(\Z,X)$.  The relative 
groups $\S_*(X)$ in the algebraic surgery exact sequence of $X$
$$\dots \to H_n(X;\L_{\bullet}) \xrightarrow[]{A} L_n(\Z[\pi_1(X)]) \to
\S_n(X) \to H_{n-1}(X;\L_{\bullet}) \to \dots$$
are the cobordism groups of quadratic Poincar\'e complexes over $(\Z,X)$
which assemble to contractible quadratic Poincar\'e complexes over
$\Z[\pi_1(X)]$.

The total surgery obstruction $s(X) \in \S_n(X)$ of an $n$-dimensional
simple Poincar\'e complex $X$ is the cobordism class of a quadratic
Poincar\'e complex over $(\Z,X)$ with contractible assembly
over $\Z[\pi_1(X)]$, which measures the homotopy invariant
part of the failure of the link of each simplex in $X$ to be a homology
sphere.  The total surgery obstruction is $s(X)=0$ if
(and for $n \geqslant 5$ only if) $X$ is simple homotopy equivalent to an
$n$-dimensional topological manifold.

The Browder-Novikov-Sullivan-Wall surgery exact sequence for an
$n$-dimensional topological manifold $M$ with $n \geqslant 5$
$$\dots \to L_{n+1}(\Z[\pi_1(M)]) \to \S^{TOP}(M) \to [M,G/TOP] \to
L_n(\Z[\pi_1(M)])$$
is identified with the corresponding portion of the algebraic surgery
exact sequence
$$\dots \to L_{n+1}(\Z[\pi_1(M)]) \to \S_{n+1}(M) \to 
H_n(M;\L_{\bullet}) \xrightarrow[]{A} L_n(\Z[\pi_1(M)])~.$$
The structure invariant $s(h) \in \S^{TOP}(M)=\S_{n+1}(M)$ of a simple
homotopy equivalence of $n$-dimensional topological manifolds $h:N \to
M$ is the cobordism class of an $n$-dimensional quadratic 
Poincar\'e complex in $(\Z,M)$ with contractible assembly over
$\Z[\pi_1(M)]$, which measures the homotopy invariant part of the
failure of the point inverses $h^{-1}(x)$ ($x \in M$) to be acyclic. 
The structure invariant is $s(h)=0$ if (and for $n
\geqslant 5$ only if) $h$ is homotopic to a homeomorphism.
\end{abstract}

\vfill

{\it Keywords:} surgery exact sequence, structure set, total surgery
obstruction

{\it AMS numbers:} 57R67, 57P10, 57N65

\newpage
\thispagestyle{empty}
\tableofcontents

\newpage
\setcounter{page}{1}

\section{Introduction}

The structure set of a differentiable $n$-dimensional manifold $M$ is
the set $\S^O(M)$ of equivalence classes of pairs $(N,h)$ with $N$ a
differentiable manifold and $h:N \to M$ a simple homotopy equivalence, 
subject to
$(N,h) \sim (N',h')$ if there exist a diffeomorphism $f:N \to N'$ and a
homotopy $f \simeq h'f:N \to M$.  The differentiable structure set was
first computed for $N=S^n$ ($n \geqslant 5$), with $\S^O(S^n)=\Theta^n$ the
Kervaire-Milnor group of exotic spheres.  In this case the structure
set is an abelian group, since the connected sum of homotopy
equivalences $h_1:N_1 \to S^n$, $h_2:N_2 \to S^n$ is a homotopy
equivalence
$$h_1 \# h_2~:~N_1 \# N_2 \to S^n \# S^n ~=~S^n~.$$
\indent The Browder-Novikov-Sullivan-Wall theory for the classification
of manifold structures within the simple homotopy type of an
$n$-dimensional differentiable manifold $M$ with $n \geqslant 5$
fits $\S^O(M)$ into an exact sequence of pointed sets
$$\dots \to L_{n+1}(\Z[\pi_1(M)]) \to \S^O(M) \to [M,G/O] \to
L_n(\Z[\pi_1(M)])$$
corresponding to the two stages of the obstruction theory for deciding
if a simple homotopy equivalence $h:N \to M$ is homotopic to a
diffeomorphism:
\begin{itemize}
\item[(i)] The primary obstruction in $[M,G/O]$ to the extension of $h$
to a normal bordism $(f,b):(W;M,N) \to M \times ([0,1];\{0\},\{1\})$
with $f\vert=1:M \to M$. Here $G/O$ is the classifying space for
fibre homotopy trivialized vector bundles, and $[M,G/O]$ is 
identified with the bordism of normal maps $M' \to M$ by the 
Browder-Novikov transversality construction.
\item[(ii)] The secondary obstruction $\sigma_*(f,b) \in
L_{n+1}(\Z[\pi_1(M)])$ to performing surgery on $(f,b)$ to make
$(f,b)$ a simple homotopy equivalence, which depends on the choice of
solution in (i).  Here, it is necessary to use the version of the
$L$-groups $L_*(\Z[\pi_1(X)])$ in which modules are based and
isomorphisms are simple, in order to take advantage of the
$s$-cobordism theorem.
\end{itemize}
The Whitney sum of vector bundles makes $G/O$ an $H$-space
(in fact an infinite loop space), so that
$[M,G/O]$ is an abelian group. However, the surgery obstruction function
$[M,G/O] \to L_n(\Z[\pi_1(M)])$ is not a morphism of groups, and
in general the differentiable structure set $\S^O(M)$ does not have
a group structure (or at least is not known to have), abelian or otherwise.

The 1960's development of surgery theory culminated in the work of
Kirby and Siebenmann \cite{KirbySiebenmann} on high-dimensional
topological manifolds, which revealed both a striking similarity and a
striking difference between the differentiable and topological
catgeories.  Define the structure set of a topological $n$-dimensional
manifold $M$ exactly as before, to be the set $\S^{TOP}(M)$ of
equivalence classes of pairs $(N,h)$ with $N$ a topological manifold
and $h:N \to M$ a simple homotopy equivalence, subject to $(N,h) \sim
(N',h')$ if there exist a homeomorphism $f:N \to N'$.  The similarity
is that again there is a surgery exact sequence for $n \geqslant 5$
$$\dots \to L_{n+1}(\Z[\pi_1(M)]) \to \S^{TOP}(M) \to [M,G/TOP] \to
L_n(\Z[\pi_1(M)]) \eqno{(*)}$$
corresponding to a two-stage obstruction theory for deciding if a simple
homotopy equivalence is homotopic to a homeomorphism, with $G/TOP$ the
classifying space for fibre homotopy trivialized topological block
bundles.  The difference is that the topological structure set
$\S^{TOP}(M)$ has an abelian group structure and $G/TOP$ has an
infinite loop space structure with respect to which $(*)$ is an exact
sequence of abelian groups.  Another difference is given by the
computation $\S^{TOP}(S^n)=0$, which is just a restatement of the
generalized Poincar\'e conjecture in the topological category :
for $n \geqslant 5$ every homotopy equivalence $h:M^n \to S^n$ is homotopic to
a homeomorphism.

Originally, the abelian group structure on $(*)$ was suggested by
the characteristic variety theorem of Sullivan \cite{Sullivan}
on the homotopy type of $G/TOP$, including the computation
$\pi_*(G/TOP)=L_*(\Z)$.  Next, Quinn \cite{Quinn1} proposed 
that the surgery obstruction function should be factored as the
composite
$$[M,G/TOP]~=~H^0(M;\underline{G/TOP})~\cong~H_n(M;\underline{G/TOP})
\xrightarrow[]{A} L_n(\Z[\pi_1(M)])$$
with $\underline{G/TOP}$ the simply-connected surgery spectrum with 0th
space $G/TOP$, identifying the topological structure sequence with the
homotopy exact sequence of a geometrically defined spectrum-level
assembly map.  This was all done in Ranicki \cite{Ranicki1979},
\cite{Ranicki1992}, but with algebra instead of geometry.

The algebraic theory of surgery was used in \cite{Ranicki1992}
to define the algebraic surgery exact sequence of abelian groups 
for any space $X$
$$\dots \to H_n(X;\L_{\bullet}) \xrightarrow[]{A} L_n(\Z[\pi_1(X)]) \to
\S_n(X) \to H_{n-1}(X;\L_{\bullet}) \to \dots~.\eqno{(**)}$$
The expression of the $L$-groups $L_*(\Z[\pi_1(X)])$ as the cobordism
groups of quadratic Poincar\'e complexes over $\Z[\pi_1(X)]$ (recalled
in the notes on the foundations of algebraic surgery) was extended to
$H_*(X;\L_{\bullet})$ and $\S_*(X)$, using quadratic Poincar\'e
complexes in categories containing much more of the topology of $X$
than just the fundamental group $\pi_1(X)$.  The topological surgery
exact sequence of an $n$-dimensional manifold $M$ with $n \geqslant 5$ was
shown to be in bijective correspondence with the corresponding portion
of the algebraic surgery sequence, including an explicit bijection
$$s~:~\S^{TOP}(M) \to \S_{n+1}(M)~;~h \mapsto s(h)~.$$
The structure invariant $s(h) \in \S_{n+1}(M)$ of a simple homotopy
equivalance $h:N \to M$ measures the chain level cobordism failure of
the point inverses $h^{-1}(x) \subset N$ ($x \in M$) to be points.

The Browder-Novikov-Sulivan-Wall surgery theory deals both with the
existence and uniqueness of manifolds in the simple homotopy type of a
geometric simple $n$-dimensional Poincar\'e complex $X$ with $n \geq
5$.  Again, this was first done for differentiable manifolds, and then
extended to topological manifolds, with a two-stage obstruction :
\begin{itemize}
\item[(i)] The primary obstruction in $[X,B(G/TOP)]$ to the existence
of a normal map $(f,b):M \to X$, with $b:\nu_M \to \widetilde{\nu}_X$ a
bundle map from the stable normal bundle $\nu_M$ of $M$ to a topological
reduction $\widetilde{\nu}_X:X \to BTOP$ of the Spivak normal fibration
$\nu_X:X \to BG$.
\item[(ii)] The secondary obstruction $\sigma_*(f,b) \in
L_n(\Z[\pi_1(M)])$ to performing surgery to make $(f,b)$ a simple
homotopy equivalence, which depends on the choice of solution in (i). 
\end{itemize}
For $n \geqslant 5$ $X$ is simple homotopy equivalent to a topological manifold
if and only if there exists a topological reduction $\widetilde{\nu}_X:X \to BTOP$ 
for which the corresponding normal map $(f,b):M \to X$ has surgery
obstruction $\sigma_*(f,b)=0$.
The algebraic surgery exact sequence $(**)$ unites the two stages into a
single invariant, the total surgery obstruction 
$s(X) \in \S_n(X)$, which measures the chain level cobordism
failure of the points $x \in X$ to have Euclidean neighbourhoods.
For $\pi_1(X)=\{1\}$, $n=4k$
the condition $s(X)=0\in \S_{4k}(X)$ is precisely the Browder condition
that there exist a topological reduction $\widetilde{\nu}_X$ for which the 
signature of $X$ is given by the Hirzebruch formula
$${\rm signature}(X)~=~\langle {\mathcal L}(-\widetilde{\nu}_X),[X] \rangle
\in \Z~.$$

\section{Geometric Poincar\'e assembly}

This section describes the assembly for geometric Poincar\'e
bordism, setting the scene for the use of quadratic Poincar\'e 
bordism in the assembly map in algebraic $L$-theory.
In both cases assembly is the passage from objects with local
Poincar\'e duality to objects with global Poincar\'e duality.

Given a space $X$ let $\Omega^P_n(X)$ be the bordism group of maps
$f:Q \to X$ from $n$-dimensional geometric Poincar\'e complexes $Q$.
The functor $X \mapsto \Omega^P_*(X)$ is homotopy invariant.
If $X=X_1\cup_YX_2$ it is not in general possible to make $f:Q \to X$
Poincar\'e transverse at $Y \subset X$, i.e. $f^{-1}(Y) \subset Q$
will not be an $(n-1)$-dimensional geometric Poincar\'e complex.
Thus $X \mapsto \Omega^P_*(X)$ does not have Mayer-Vietoris sequences,
and is not a generalized homology theory. 
The general theory
of Weiss and Williams \cite{WeissWilliams} provides a generalized
homology theory $X \mapsto H_*(X;\Omega^P_{\bullet})$ with an assembly map
$A:H_*(X;\Omega^P_{\bullet}) \to \Omega^P_*(X)$. However, it
is possible to obtain $A$ by a direct geometric construction :
$H_n(X;\Omega^P_{\bullet})$ is the bordism
group of Poincar\'e transverse maps $f:Q \to X$ from $n$-dimensional
Poincar\'e complexes $Q$, and $A$ forgets the transversality.
The coefficient spectrum $\Omega^P_{\bullet}$ is such that
$$\pi_*(\Omega^P_{\bullet})~=~\Omega^P_*(\{{\rm pt.}\})~,$$ 
and may be constructed using geometric Poincar\'e $n$-ads.

In order to give a precise geometric description of 
$H_n(X;\Omega^P_{\bullet})$ it is convenient to assume that $X$ is
the polyhedron of a finite simplicial complex (also denoted $X$).
The {\it dual cell} of a simplex $\sigma \in X$ is the subcomplex of
the barycentric subdivision $X'$
$$D(\sigma,X)~=~\{\widehat{\sigma}_0\widehat{\sigma}_1 \dots
\widehat{\sigma}_n \,\vert\, \sigma \leqslant \sigma_0 < \sigma_1 < \dots <
\sigma_n\} \subset X'~,$$
with boundary the subcomplex
$$\partial D(\sigma,X)~=~\bigcup_{\tau>\sigma}D(\tau,X)~=~
\{\widehat{\tau}_0\widehat{\tau}_1 \dots \widehat{\tau}_n \,\vert\,
\sigma < \tau_0 < \tau_1 < \dots < \tau_n\} \subset D(\sigma,X)~.$$
Every map $f:M \to X$ from an $n$-manifold $M$ can be made transverse
across the dual cells, meaning that for each $\sigma \in X$
$$(M(\sigma),\partial M(\sigma))~=~f^{-1}(D(\sigma,X),\partial D(\sigma,X))$$
is an $(n-\vert \sigma \vert)$-dimensional manifold with boundary.
Better still, for an $n$-dimensional $PL$ manifold $M$ every simplicial map 
$f:M \to X'$ is already transverse in this sense, by a result of Marshall
Cohen.

A map $f:Q \to X$ is {\it $n$-dimensional Poincar\'e transverse} if for
each $\sigma\in X$
$$(Q(\sigma),\partial Q(\sigma))~=~f^{-1}(D(\sigma,X),\partial D(\sigma,X))$$
is an $(n-\vert \sigma \vert)$-dimensional geometric Poincar\'e pair.

\noindent{\it Proposition.} $H_n(X;\Omega^P_{\bullet})$ is the bordism
group of Poincar\'e transverse maps $Q \to X$ from $n$-dimensional
geometric Poincar\'e complexes.\hfill\qed

It is worth noting that 
\begin{itemize}
\item[(i)] The identity $1:X \to X$ is $n$-dimensional Poincar\'e
transverse if and only if $X$ is an $n$-dimensional homology manifold.
\item[(ii)] If a map $f:Q \to X$ is $n$-dimensional Poincar\'e transverse
then $Q$ is an $n$-dimensional geometric Poincar\'e complex. The global
Poincar\'e duality of $Q$ is assembled
from the local Poincar\'e dualities of $(Q(\sigma),\partial Q(\sigma))$.
For $f=1:Q=X \to X$ this is the essence of 
Poincar\'e's original proof of his duality for a homology manifold.
\end{itemize}

The {\it Poincar\'e structure group} $\S^P_n(X)$ is the relative group
in the {\it geometric Poincar\'e surgery exact sequence}
$$\dots \to H_n(X;\Omega_{\bullet}^P) \xrightarrow[]{A} \Omega^P_n(X) \to \S^P_n(X) \to
H_{n-1}(X;\Omega_{\bullet}^P) \to \dots~,$$
which is the cobordism group of maps $(f,\partial f):(Q,\partial Q) \to X$
from $n$-dimensional  Poincar\'e pairs $(Q,\partial Q)$ with
$\partial f: \partial Q \to X$ Poincar\'e transverse. 
The {\it total Poincar\'e surgery obstruction} of an $n$-dimensional geometric 
Poincar\'e complex $X$ is the image $s^P(X) \in \S^P_n(X)$ of 
$(1:X \to X) \in \Omega_n^P(X)$, with $s^P(X)=0$ if and only if
there exists an $\Omega_{\bullet}^P$-coefficient fundamental class
$[X]_P \in H_n(X;\Omega_{\bullet}^P)$ with $A([X]_P)=(1:X \to X) \in 
\Omega^P_n(X)$.

In fact, it follows from the Levitt-Jones-Quinn-Hausmann-Vogel
Poincar\'e bordism theory that $\S^P_n(X)=\S_n(X)$ for $n \geqslant 5$,
and that $s^P(X)=0$ if and only if $X$ is homotopy equivalent to an
$n$-dimensional topological manifold.  The geometric Poincar\'e bordism
approach to the structure sets and total surgery obstruction is
intuitive, and has the virtue(?) of dispensing with the algebra
altogether.  Maybe it even applies in the low dimensions $n=3,4$. 
However, at present our understanding of the Poincar\'e bordism theory
is not good enough to use it for foundational purposes.  So back to the
algebra!

\section{The algebraic surgery exact sequence}

This section constructs the quadratic $L$-theory assembly map $A$
and the algebraic surgery exact sequence
$$\dots \to H_n(X;\L_{\bullet}) \xrightarrow[]{A} L_n(\Z[\pi_1(X)]) \to \S_n(X) \to
H_{n-1}(X;\L_{\bullet}) \to \dots\eqno{(**)}$$
for a finite simplicial complex $X$. A `$(\Z,X)$-module' is a based f.g.
free $\Z$-module in which every basis element is associated to a simplex
of $X$. The construction of $(**)$ makes use
of a chain complex duality on the $(\Z,X)$-module category $\A(\Z,X)$.

The quadratic $L$-spectrum $\L_{\bullet}$ is 
1-connective, with connected 0th space $\L_0\simeq G/TOP$ and
$$\pi_n(\L_{\bullet})~=~\pi_n(\L_0)~=~L_n(\Z)~=~\begin{cases} 
\Z&\text{if $n\equiv 0(\bmod\,4)$ (signature)/8}\\
\Z_2&\text{if $n\equiv 2(\bmod\,4)$ (Arf invariant)}\\
0&\text{otherwise}
\end{cases}$$
for $n \geqslant 1$. From the algebraic point of view it is easier to start 
with the 0-connective quadratic $L$-spectrum 
$\overline{\L}_{\bullet}=\L_{\bullet}(\Z)$, such that
$$\pi_n(\overline{\L}_{\bullet})~=~\begin{cases} L_n(\Z)&\hbox{if
$n\geqslant 0$} \cr
0&\hbox{if $n \leqslant -1$}
\end{cases}$$
with disconnected 0th space $\overline{\L}_0 \simeq L_0(\Z) \times G/TOP$.
The two
spectra are related by a fibration sequence
$\L_{\bullet} \to \overline{\L}_{\bullet} \to \K(L_0(\Z))$
with $\K(L_0(\Z))$ the Eilenberg-MacLane spectrum of $L_0(\Z)$. 

The algebraic surgery exact sequence was constructed in Ranicki
\cite{Ranicki1992} using the $(\Z,X)$-module category of Ranicki and
Weiss \cite{RanickiWeiss}.  (This is a rudimentary version of
controlled topology, cf.  Ranicki \cite{Ranicki1999}).

A {\it $(\Z,X)$-module} is a direct sum of based f.g. free $\Z$-modules
$$B~=~\sum\limits_{\sigma \in X}B(\sigma)~.$$
A $(\Z,X)$-module morphism $f:B \to C$ is a $\Z$-module morphism such that
$$f(B(\sigma)) \subseteq\sum\limits_{\tau \geqslant \sigma}C(\tau)~,$$
so that the matrix of $f$ is upper triangular. 
A $(\Z,X)$-module chain map $f:B \to C$ is a chain equivalence if and
only if each $f(\sigma,\sigma):B(\sigma) \to C(\sigma)$ $(\sigma \in X)$
is a $\Z$-module chain equivalence. 
The universal covering projection $p:\widetilde{X} \to X$ is used 
to define the $(\Z,X)$-module assembly functor
$$A~:~\A(\Z,X) \to \A(\Z[\pi_1(X)])~;~ B \mapsto 
\sum\limits_{\widetilde{\sigma} \in \widetilde{X}}B(p\widetilde{\sigma})$$
with $\A(\Z,X)$ the category of $(\Z,X)$-modules and
$\A(\Z[\pi_1(X)])$ the category of based f.g. free $\Z[\pi_1(X)]$-modules.
In the language of sheaf theory $A=q_!p^!$ (cf. Verdier \cite{Verdier}),
with $q:\widetilde{X} \to \{\text{pt.}\}$.

The involution $g \mapsto \overline{g}=g^{-1}$ on $\Z[\pi_1(X)]$
extends in the usual way to a duality involution on $\A(\Z[\pi_1(X)])$,
sending a based f.g. free $\Z[\pi_1(X)]$-module $F$ to the dual f.g.
free $\Z[\pi_1(X)]$-module $F^*={\rm Hom}_{\Z[\pi_1(X)]}(F,\Z[\pi_1(X)])$.
Unfortunately, it is not possible to define a duality involution on
$\A(\Z,X)$ (since the transpose of an upper triangular matrix is a
lower triangular matrix). See Chapter 5 of Ranicki \cite{Ranicki1992}
for the construction of a `chain duality' on $\A(\Z,X)$ and
of the $L$-groups $L_*(\A(\Z,X))$. The chain duality associates 
to a chain complex $C$ in $\A(\Z,X)$ a chain complex $TC$ in $\A(\Z,X)$ with
$$TC(\sigma)_r~=~\sum\limits_{\tau \geqslant \sigma}
\text{\rm Hom}_{\Z}(C_{-\vert \sigma \vert -r}(\tau),\Z)~.$$
\noindent{\it Example.} The simplicial chain complex $C(X')$ is a
$(\Z,X)$-module chain complex, with assembly $A(C(X'))$ $\Z[\pi_1(X)]$-module
chain equivalent to $C(\widetilde{X})$. The chain dual $TC(X')$ is
$(\Z,X)$-module chain equivalent to the simplicial cochain complex 
$D={\rm Hom}_{\Z}(C(X),\Z)^{-*}$, with assembly $A(D)$ 
which is $\Z[\pi_1(X)]$-module
chain equivalent to $C(\widetilde{X})^{-*}$.\hfill\qed

The quadratic $L$-group $L_n(\A(\Z,X))$ is the cobordism group of
$n$-dimensional quadratic Poincar\'e complexes $(C,\psi)$ in $\A(\Z,X)$.

\noindent{\it Proposition.} (\cite{Ranicki1992}, 14.5)
The functor $X \mapsto L_*(\A(\Z,X))$ is the generalized homology theory 
with $\L_{\bullet}(\Z)$-coefficients
$$L_*(A(\Z,X))~=~H_*(X;\L_{\bullet}(\Z))~.\eqno{\hbox{\qed}}$$

The coefficient spectrum $\overline{\L}_{\bullet}=\L_{\bullet}(\Z)$ is
the special case $R=\Z$ of a general construction.  For any ring with
involution $R$ there is a 0-connective spectrum $\L_{\bullet}(R)$ such
that
$$\pi_*(\overline{\L}_{\bullet}(R))~=~L_*(R)~,$$
which may be constructed using  quadratic Poincar\'e $n$-ads over $R$.

The assembly functor $A:\A(\Z,X) \to \A(\Z[\pi_1(X)])$
induces assembly maps in the quadratic $L$-groups, which fit into
the {\it 4-periodic algebraic surgery exact sequence}
$$\dots \to H_n(X;\overline{\L}_{\bullet}) \xrightarrow[]{A} L_n(\Z[\pi_1(X)])
\to \overline{\S}_n(X) \to H_{n-1}(X;\overline{\L}_{\bullet}) \to \dots $$
with the {\it 4-periodic algebraic structure set} $\overline{\S}_n(X)$
the cobordism group of $(n-1)$-dimensional quadratic Poincar\'e
complexes $(C,\psi)$ in $\A(\Z,X)$ such that the assembly $A(C)$ is a
simple contractible based f.g.  free $\Z[\pi_1(X)]$-module chain
complex.  (See section 4.5 for the geometric interpretation).  A
priori, an element of the relative group $\overline{\S}_n(X)=\pi_n(A)$
is an $n$-dimensional quadratic $\Z[\pi_1(X)]$-Poincar\'e pair $(C \to
D,(\delta\psi,\psi))$ in $\A(\Z,X)$.  Using this as data for algebraic
surgery results in an $(n-1)$-dimensional quadratic Poincar\'e complex
$(C',\psi')$ in $\A(\Z,X)$ such that the assembly $A(C')$ is a simple
contractible based f.g.  free $\Z[\pi_1(X)]$-module chain complex.

Killing $\pi_0(\overline{\L}_{\bullet})=L_0(\Z)$ in 
$\overline{\L}_{\bullet}$ results in the 1-connective
spectrum $\L_{\bullet}$, and the {\it algebraic surgery exact sequence}
$$\dots \to H_n(X;\L_{\bullet}) \xrightarrow[]{A} L_n(\Z[\pi_1(X)])
\to \S_n(X) \to H_{n-1}(X;\L_{\bullet}) \to \dots\eqno{(**)}$$
with $\S_n(X)$ the {\it algebraic structure set}. The two types
of structure set are related by an exact sequence
$$\dots \to H_n(X;L_0(\Z)) \to \S_n(X) \to \overline{\S}_n(X) \to
H_{n-1}(X;L_0(\Z))\to \dots~.$$

\section{The structure set and the total surgery obstruction}

This chapter states the results in Chapters 16,17,18
of Ranicki \cite{Ranicki1992} on the $L$-theory orientation
of topology, the total surgery obstruction and the structure set.

The algebraic theory of surgery fits the homotopy category
of topological manifolds of dimension $\geqslant 5$
into a pullback square 
$$\xymatrix@R+10pt{
\{\text{topological manifolds}\} \ar[r] \ar[d] &
\{\text{local algebraic Poincar\'e complexes}\} \ar[d] \\
\{\text{geometric Poincar\'e complexes}\} \ar[r] & 
\{\text{global algebraic Poincar\'e complexes}\}
}$$
where local means $\A(\Z,X)$ and global means $\A(\Z[\pi_1(X)])$.
In words~: the homotopy type of a topological manifold is 
the homotopy type of a geometric Poincar\'e complex with a local
algebraic Poincar\'e structure.

\subsection{The $L$-theory orientation of topological block bundles}

The topological $k$-block bundles of Rourke and Sanderson
\cite{RourkeSanderson} are topological analogues of vector bundles.  By
analogy with the classifying spaces $BO(k)$, $BO$ for vector bundles
there are classifying spaces $B\widetilde{TOP}(k)$ for topological
block bundles, and a stable classifying space $BTOP$.  It is known from
the work of Sullivan \cite{Sullivan} and Kirby-Siebenmann
\cite{KirbySiebenmann} that the classifying space for fibre homotopy
trivialized topological block bundles
$$G/TOP~=~\hbox{\rm homotopy fibre}(BTOP \to BG)$$
has homotopy groups $\pi_*(G/TOP)=L_*(\Z)$.
A map $S^n \to G/TOP$ classifies a topological block bundle
$\eta:S^n \to B\widetilde{TOP}(k)$ with a fibre homotopy trivialization
$J\eta \simeq \{*\}:S^n \to BG(k)$ ($k \geqslant 3$). 
The isomorphism $\pi_n(G/TOP) \to L_n(\Z)$ is defined by sending 
$S^n \to G/TOP$ to the surgery obstruction $\sigma_*(f,b)$ of
the corresponding normal map $(f,b):M \to S^n$ from a topological
$n$-dimensional manifold $M$, with $b:\nu_M \to \nu_{S^n}\oplus\eta$.
Sullivan \cite{Sullivan} proved that $G/TOP$ and $BO$ have the
same homotopy type localized away from 2 
$$G/TOP[1/2]~\simeq~BO[1/2]~.$$
(The localization $\Z[1/2]$ is the subring 
$\{\ell/2^m\,\vert\,\ell \in \Z,m \geqslant 0\} \subset \Q$ 
obtained from $\Z$ by inverting 2.  The
localization $X[1/2]$ of a space $X$ is a space such that
$$\pi_*(X[1/2])~=~\pi_*(X)\otimes_{\Z}\Z[1/2]~.$$
Thus $X \mapsto X[1/2]$ kills all the 2-primary torsion in $\pi_*(X)$.)

Let $\L^{\bullet}=\L(\Z)^{\bullet}$ be the symmetric $L$-spectrum of $\Z$,
with homotopy groups
$$\pi_n(\L^{\bullet})~=~L^n(\Z)~=~\begin{cases} 
\Z&\text{if $n\equiv 0(\bmod\,4)$ (signature)}\\
\Z_2&\text{if $n\equiv 1(\bmod\,4)$ (deRham invariant)}\\
0&\text{otherwise}~.
\end{cases}$$
The hyperquadratic $\L$-spectrum of $\Z$ is defined by
$$\widehat{\L}^{\bullet}~=~{\rm cofibre}(1+T:\L_{\bullet} \to \L^{\bullet})~.$$
It is 0-connective, fits into a (co)fibration sequence of spectra
$$\dots \to \L_{\bullet} \xrightarrow[]{1+T} \L^{\bullet} \to
\widehat{\L}^{\bullet} \to \Sigma \L_{\bullet} \to \dots~,$$
and has homotopy groups
$$\pi_n(\widehat{\L}^{\bullet})~=~\begin{cases} 
\Z&\hbox{if $n=0$}\\
\Z_8&\hbox{if $n \equiv 0(\bmod\,4)$ and $n>0$}\\
\Z_2&\hbox{if $n \equiv 2,3(\bmod\,4)$}\\
0& \hbox{if $n \equiv 1(\bmod\,4)$~.}
\end{cases}$$

An {\it $h$-orientation} of a spherical fibration $\nu:X \to BG(k)$ with
respect to a ring spectrum $h$ is an $h$-coefficient Thom class
in the reduced $h$-cohomology $U \in \dot h^k(T(\nu))$ of the Thom space
$T(\nu)$, i.e. a $\dot h$-cohomology class 
which restricts to $1 \in \dot h^k(S^k)=\pi_0(h)$ over each $x \in X$.
\smallskip

\noindent{\it Theorem} (\cite{Ranicki1992}, 16.1) 
(i) The $0$th space $\L_0$ of $\L_{\bullet}$ is homotopy equivalent to $G/TOP$
$$\L_0~\simeq~G/TOP~.$$
(ii) Every topological $k$-block bundle $\nu:X \to B\widetilde{TOP}(k)$
has a canonical $\L^{\bullet}$-orientation
$$U_{\nu} \in \dot{H}^k(T(\nu);\L^{\bullet})~.$$
(iii) Every $(k-1)$-spherical fibration $\nu:X \to BG(k)$ has a 
canonical $\widehat{\L}^{\bullet}$-orientation 
$$\widehat{U}_{\nu} \in \dot H^k(T(\nu);\widehat{\L}^{\bullet})~,$$
with $\dot H$ denoting reduced cohomology.
The {\it topological reducibility obstruction}
$$t(\nu)~=~\delta(\widehat{U}_{\nu}) \in \dot{H}^{k+1}(T(\nu);\L_{\bullet})$$ 
is such that $t(\nu)=0$ if and only if $\nu$ admits a  topological block bundle
reduction $\widetilde{\nu}:X \to B\widetilde{TOP}(k)$. Here, $\delta$
is the connecting map in the exact sequence
$$\dots \to \dot H^k(T(\nu);\L_{\bullet}) \to 
\dot H^k(T(\nu);\L^{\bullet}) \to \dot H^k(T(\nu);\widehat{\L}^{\bullet}) 
\xrightarrow[]{\delta} \dot H^{k+1}(T(\nu);\L_{\bullet}) \to \dots~.$$
The topological block bundle reductions of $\nu$ are in one-one correspondence
with lifts of $\widehat{U}_{\nu}$ to a $\L^{\bullet}$-orientation 
$U_{\nu} \in \do H^k(T(\nu);\L^{\bullet})$.\hfill\qed

\noindent{\it Example.} Rationally, the symmetric $L$-theory orientation 
of $\nu:X \to B\widetilde{TOP}(k)$ is the ${\mathcal L}$-genus
$$U_{\nu}\otimes \Q~=~{\mathcal L}(\nu) \in \dot{H}^k(T(\nu);\L^{\bullet})\otimes \Q~=~
H^{4*}(X;\Q)~.\eqno{\hbox{\qed}}$$

\noindent{\it Example.} Localized away from 2, the symmetric $L$-theory
orientation of $\nu:X \to B\widetilde{TOP}(k)$
is the $KO[1/2]$-orientation of Sullivan \cite{Sullivan}
$$U_{\nu}[1/2]~=~\Delta_\nu \in \dot{H}^k(T(\nu);\L^{\bullet})[1/2]~=~
\widetilde{KO}^k(T(\nu))[1/2]~.\eqno{\hbox{\qed}}$$

\subsection{The total surgery obstruction}

The {\it total surgery obstruction} $s(X) \in \S_n(X)$ of an
$n$-dimensional geometric Poincar\'e complex $X$ is the cobordism class
of the $\Z[\pi_1(X)]$-contractible $(n-1)$-dimensional quadratic
Poincar\'e complex $(C,\psi)$ in $\A(\Z,X)$ with
$C=C([X]\cap-:C(X)^{n-*} \to C(X'))_{*+1}$, using the dual cells in the
barycentric subdivision $X'$ to regard the simplicial chain complex
$C(X')$ as a chain complex in $\A(\Z,X)$.
\medskip

\noindent{\it Theorem} (\cite{Ranicki1992}, 17.4) The total surgery
obstruction is such that $s(X)=0 \in \S_n(X)$ if (and for $n \geqslant 5$
only if) $X$ is homotopy equivalent to an $n$-dimensional topological
manifold.\\
\noindent{\it Proof} A regular neighbourhood $(W,\partial W)$ of 
an embedding $X \subset S^{n+k}$ ($k$ large) gives a Spivak normal fibration
$$S^{k-1} \to \partial W \to W~\simeq~X$$
with Thom space $T(\nu)=W/\partial W$ $S$-dual to $X_+=X\cup \{{\rm pt.}\}$.
The total surgery obstruction $s(X) \in \S_n(X)$ has image the
topological reducibility obstruction
$$t(\nu) \in \dot H^{k+1}(T(\nu);\L_{\bullet})~\cong~ 
H_{n-1}(X;\L_{\bullet})~.$$
Thus $s(X)$ has image $t(\nu)=0 \in H_{n-1}(X;\L_{\bullet})$
if and only if $\nu$ admits a topological block bundle
reduction $\widetilde{\nu}:X \to B\widetilde{TOP}(k)$,
in which case the topological version of the Browder-Novikov
transversality construction applied to the degree 1 map
$\rho:S^{n+k} \to T(\nu)$ gives a normal map 
$(f,b)=\rho\vert:M=f^{-1}(X) \to X$. The surgery obstruction
$\sigma_*(f,b) \in L_n(\Z[\pi_1(X)])$ has image 
$$[\sigma_*(f,b)]~=~s(X) \in
{\rm im}(L_n(\Z[\pi_1(X)]) \to \S_n(X))~=~
{\rm ker}(\S_n(X) \to H_{n-1}(X;\L_{\bullet}))$$
The total surgery obstruction is $s(X)=0$ if and only if
there exists a reduction $\widetilde{\nu}$ with $\sigma_*(f,b)=0$.
\hfill\qed

\noindent{\it Example.} 
For a simply-connected space $X$ the assembly map 
$A:H_*(X;\L_{\bullet}) \to L_*(\Z)$ is onto, so that
$$\S_n(X)~=~{\rm ker}(A:H_{n-1}(X;\L_{\bullet}) \to L_{n-1}(\Z))~
=~\dot H_{n-1}(X;\L_{\bullet})~ ,$$
with $\dot H$ denoting reduced homology.
The total surgery obstruction $s(X) \in \S_n(X)$ 
of a simply-connected $n$-dimensional
geometric Poincar\'e complex $X$ is just the obstruction to
the topological reducibility of the Spivak normal fibration
$\nu_X:X \to BG$.\hfill\qed

There are also relative and rel $\partial$ versions of the total
surgery obstruction.

For any pair of spaces $(X,Y \subseteq X)$ let $\S_n(X,Y)$ be the
relative groups in the exact sequence
$$\dots \to H_n(X,Y;\L_{\bullet}) \xrightarrow[]{A} L_n(\Z[\pi_1(Y)]
\to \Z[\pi_1(X)])
\to \S_n(X,Y) \to H_{n-1}(X,Y;\L_{\bullet}) \to \dots~.$$
The {\it relative total surgery obstruction} $s(X,Y) \in \S_n(X,Y)$ of an 
$n$-dimensional geometric Poincar\'e pair 
is such that $s(X,Y)=0$ if (and for $n \geqslant 6$ only if) $(X,Y)$
is homotopy equivalent to an $n$-dimensional topological
manifold with boundary $(M,\partial M)$.
In the special case $\pi_1(X)=\pi_1(Y)$
$$s(X,Y) \in \S_n(X,Y) = H_{n-1}(X,Y;\L_{\bullet})$$ 
is just the obstruction to the topological reducibility of the Spivak
normal fibration $\nu_X:X \to BG$, which is the $\pi$-$\pi$ theorem
of Chapter 4 of Wall \cite{Wall}.

The {\it rel $\partial$ total surgery obstruction} $s_{\partial}(X,Y)
\in \S_n(X)$ of an $n$-dimensional geometric Poincar\'e pair $(X,Y)$
with manifold boundary $Y$ is such that $s_{\partial}(X,Y)=0$ if (and
for $n \geqslant 5$ only if) $(X,Y)$ is homotopy equivalent rel $\partial$
to an $n$-dimensional manifold with boundary.

\subsection{The $L$-theory orientation of topological manifolds}

An $n$-dimensional geometric Poincar\'e complex $X$ determines
a symmetric $\Z[\pi_1(X)]$-Poincar\'e complex $(C(X'),\phi)$ in $\A(\Z,X)$,
with assembly the usual symmetric Poinacar\'e complex
$(C(\widetilde{X}),\phi(\widetilde{X}))$
representing the symmetric signature $\sigma^*(X) \in L^n(\Z[\pi_1(X)])$.
\medskip

\noindent{\it Example.} For $n=4k$ $\sigma^*(M) \in
L^{4k}(\Z[\pi_1(M)])$ has image
$${\rm signature}(X)~=~{\rm signature}(H^{2k}(X;\Q),\cup)\in L^{4k}(\Z)=\Z~.
\eqno{\hbox{\qed}}$$

A triangulated $n$-dimensional manifold $M$ determines a symmetric 
Poincar\'e complex $(C(M'),\phi)$ in $\A(\Z,M)$.
The {\it symmetric $L$-theory orientation} of $M$ is the 
$\L^{\bullet}$-coefficient class
$$[M]_{\L}~=~(C(M'),\phi) \in L^n(\A(\Z,M))~=~H_n(M;\L^{\bullet})$$
with assembly
$$A([M]_{\L})~=~\sigma^*(M) \in L^n(\Z[\pi_1(M)])~.$$

\noindent{\it Example.} Rationally, the symmetric $L$-theory orientation
is the Poincar\'e dual of the ${\mathcal L}$-genus
$$[M]_{\L}~=~{\mathcal L}(M)\cap [M]_{\Q} \in H_n(M;\L^{\bullet})\otimes \Q~=~
H_{n-4*}(M;\Q)~=~H^{4*}(M;\Q)~.$$
Thus $A([M]_{\L})=\sigma^*(M) \in L^n(\Z[\pi_1(M)])$ is a
$\pi_1(M)$-equivariant generalization of the Hirzebruch signature theorem
for a $4k$-dimensional manifold
$${\rm signature}(M)~=~\langle {\mathcal L}(-\widetilde{\nu}_M),[M] \rangle
\in L^{4k}(\Z)~=~\Z~.\eqno{\hbox{\qed}}$$

\noindent{\it Example.} Localized away from 2, the symmetric $L$-theory
orientation is the $KO[1/2]$-orientation $\Delta(M)$ 
of Sullivan \cite{Sullivan}
$$[M]_{\L}\otimes \Z[1/2]~=~\Delta(M) \in
H_n(M;\L^{\bullet})[1/2]~=~KO_n(M)[1/2]~.\eqno{\hbox{\qed}}$$

See Chapter 16 of \cite{Ranicki1992} for the detailed definition of the
{\it visible symmetric $L$-groups} $VL^*(X)$ of a space $X$,
with the following properties :
\begin{itemize}
\item[(i)] $VL^n(X)$ is the cobordism group of $n$-dimensional
symmetric complexes $(C,\phi)$ in $\A(\Z,X)$ such that the
assembly $A(C,\phi)$ is an $n$-dimensional symmetric Poincar\'e
complex in $\A(\Z[\pi_1(X)])$, and such that each
$(C(\sigma),\phi(\sigma))$ ($\sigma \in X^{(n)}$) is a 0-dimensional
symmetric Poincar\'e complex in $\A(\Z)$.
\item[(ii)] The (covariant) functor $X \mapsto VL^*(X)$ is homotopy invariant.
\item[(iii)] The visible symmetric $L$-groups $VL^*(K(\pi,1))$ of
an Eilenberg-MacLane space $K(\pi,1)$ of a group $\pi$
are the visible symmetric $L$-groups $VL^*(\Z[\pi])$ of Weiss \cite{Weiss}.
\item[(iv)] The $VL$-groups fit into a commutative braid of exact sequences
$$\hskip5pt\xymatrix@!C@C-38pt@R-10pt{
\S_{n+1}(X)\ar[dr]\ar@/^2pc/[rr]^{}&&H_n(X;\L^{\bullet})
\ar[dr]^-{A}\ar@/^2pc/[rr]&&H_n(X;\widehat{\L}^{\bullet})\\
&H_n(X;\L_{\bullet})\ar[ur]^-{1+T} \ar[dr]^-{A}&&
VL^n(X)\ar[ur]\ar[dr]&&\\
H_{n+1}(X;\widehat{\L}^{\bullet})\ar[ur]\ar@/_2pc/[rr]&&
L_n(\Z[\pi_1(X)])\ar[ur]^-{1+T}\ar@/_2pc/[rr]_{}&& \S_n(X)}$$ 
\item[(v)] Every $n$-dimensional simple Poincar\'e complex $X$
has a visible symmetric signature $\sigma^*(X) \in VL^n(X)$ with
image the total surgery obstruction $s(X) \in \S_n(X)$.
\end{itemize}

An {\it $h$-orientation} of an $n$-dimensional Poincar\'e complex $X$
with respect to ring spectrum $h$ is an $h$-homology class
$[X]_h \in h_n(X)$ which corresponds under the $S$-duality
isomorphism $h_n(X) \cong \dot h^{k+1}(T(\nu))$ to an $h$-coefficient
Thom class $U_h \in \dot h^k(T(\nu))$ of the Spivak normal fibration
$\nu:X \to BG(k)$ ($k$ large, $X \subset S^{n+k}$).
\medskip

\noindent{\it Theorem} (\cite{Ranicki1992}, 16.16)
Every $n$-dimensional topological manifold $M$ has a canonical
$\L^{\bullet}$-orient\-ation $[M]_{\L}\in H_n(M;\L^{\bullet})$ 
with assembly
$$A([M]_{\L})~=~\sigma^*(M) \in VL^n(M)~.\eqno{\hbox{\qed}}$$

If $M$ is triangulated by a simplicial complex $K$ then
$$[M]_{\L}~=~(C,\phi)\in H_n(M;\L^{\bullet})~=~L^n(\Z,K)$$
is the cobordism class of an $n$-dimensional symmetric
Poincar\'e complex $(C,\phi)$ in $\A(\Z,K)$ with $C=C(K')$.
\medskip

\noindent{\it Example.} The canonical $\L^{\bullet}$-homology
class of an $n$-dimensional manifold $M$ is given rationally
by the Poincar\'e dual of the ${\mathcal L}(M)$-genus
${\mathcal L}(M) \in H^{4*}(M;\Q)$
$$[M]_{\L}\otimes \Q~=~{\mathcal L}(M) \cap [M]_{\Q} \in
H_n(M;\L^{\bullet})\otimes \Q~=~H_{n-4*}(M;\Q)~.
\eqno{\hbox{\qed}}$$

\noindent{\it Theorem} (\cite{Ranicki1992}, pp. 190--191) 
For $n \geqslant 5$ an $n$-dimensional simple Poincar\'e complex $X$ is
simple homotopy equivalent to an $n$-dimensional topological manifold
if and only if there exists a symmetric $L$-theory fundamental class
$[X]_{\L} \in H_n(X;\L^{\bullet})$ with assembly
$$A([X]_{\L})~=~\sigma^*(X) \in VL^n(X)~.\eqno{\hbox{\qed}}$$

In the simply-connected case $\pi_1(X)=\{1\}$ with $n=4k$
this is just :
\medskip

\noindent{\it Example.} For $k \geqslant 2$ a simply-connected
$4k$-dimensional Poincar\'e complex $X$ is homotopy equivalent to a
$4k$-dimensional topological manifold if and only if the Spivak normal
fibration $\nu_X:X \to BG$ admits a topological reduction
$\widetilde{\nu}_X:X \to BTOP$ for which the Hirzebruch signature
formula
$${\rm signature}(X)~=~\langle {\mathcal L}(-\widetilde{\nu}_X),[X] \rangle
\in L^{4k}(\Z)~=~\Z$$
holds. The if part is the topological version of the 
original result of Browder \cite{Browder} on the converse of the Hirzebruch 
signature theorem for the homotopy types of differentiable manifolds.
$$\eqno{\hbox{\qed}}$$

\subsection{The structure set}

The {\it structure invariant} of a homotopy
equivalence $h:N \to M$ of $n$-dimensional topological manifolds is
is the rel $\partial$ total surgery obstruction
$$s(h)~=~s_{\partial}(W,M \cup N) \in \S_{n+1}(W)~=~\S_{n+1}(M)$$
of the $(n+1)$-dimensional geometric Poincar\'e pair with manifold 
boundary $(W,M \cup N)$ defined by the mapping cylinder $W$ of $h$.

Here is a more direct description of the structure invariant,
in terms of the point inverses $h^{-1}(x) \subset N$ ($x \in M$).
Choose a simplicial complex $K$ with a homotopy
equivalence $g:M \to K$ such that $g$ and $gh:N \to K$ are
topologically transverse across the dual cells $D(\sigma,K) \subset K'$.
(For triangulated $M$ take $K=M$). Then $s(h)$ is the cobordism class 
$$s(h)~=~(C,\psi) \in \S_{n+1}(K)~=~\S_{n+1}(M)$$
of a $\Z[\pi_1(M)]$-contractible $n$-dimensional quadratic 
Poincar\'e complex $(C,\psi)$ in $\A(\Z,K)$ with
$$C~=~C(h:C(N) \to C(K'))_{*+1}~.$$

\noindent{\it Theorem} (\cite{Ranicki1992}, 18.3, 18.5)
(i) The structure invariant is
such that $s(h)=0 \in \S_{n+1}(M)$ if (and for $n \geqslant 5$ only if)
$h$ is homotopic to a homeomorphism.\\
(ii) The Sullivan-Wall surgery sequence of an
$n$-dimensional topological manifold $M$ with $n \geqslant 5$ is in one-one
correspondence with a portion of the algebraic surgery exact sequence,
by a bijection
$$\xymatrix@C+5pt@R+5pt{
\dots \ar[r]& L_{n+1}(\Z[\pi_1(M)]) \ar@{=}[d] \ar[r] &
\S^{TOP}(M) \ar[d]^{\displaystyle{s}}_-{\cong} \ar[r] & [M,G/TOP]
\ar[d]^{\displaystyle{t}}_-{\cong} \ar[r] & L_n(\Z[\pi_1(M)]) \ar@{=}[d]\\
\dots \ar[r]& L_{n+1}(\Z[\pi_1(M)]) \ar[r] &\S_{n+1}(M) \ar[r] &
H_n(M;\L_{\bullet}) \ar[r]^-{\displaystyle{A}}  & L_n(\Z[\pi_1(M)])}$$
The higher structure groups are the rel $\partial$ structure sets
$$\S_{n+k+1}(M)~=~\S^{TOP}_{\partial}(M \times D^k,M \times
S^{k-1})~~(k \geqslant 1)$$
of homotopy equivalences 
$(h,\partial h):(N,\partial N) \to (M \times D^k,M \times S^{k-1})$
with $\partial h : \partial N \to M \times S^{k-1}$ a homeomorphism.
\hfill\qed

\noindent{\it Example.} For a simply-connected space $M$ the assembly maps 
$A:H_*(M;\L_{\bullet}) \to L_*(\Z)$ are onto.
Thus for a simply-connected $n$-dimensional manifold $M$ 
$$\begin{array}{ll}
\S^{TOP}(M)&=~\S_{n+1}(M)\\[1ex]
&=~{\rm ker}(A:H_n(M;\L_{\bullet}) \to L_n(\Z))~=~
\dot H_n(M;\L_{\bullet})\\[1ex]
&=~{\rm ker}(\sigma_*:[M,G/TOP] \to L_n(\Z))
\end{array}$$
with $\sigma_*$ the surgery obstruction map.  The structure invariant
$s(h) \in \S^{TOP}(M)$ of a homotopy equivalence $h:N \to M$ is given
modulo 2-primary torsion by the difference of the canonical
$\L^{\bullet}$-orientations
$$s(h)[1/2]~=~(h_*[N]_{\L} - [M]_{\L},0) \in \dot
H_n(M;\L^{\bullet})[1/2]~=~\dot H_n(M;\L_{\bullet})[1/2] \oplus
H_n(M)[1/2]~.$$
Rationally, this is just the difference of the Poincar\'e duals of the
${\mathcal L}$-genera
$$\begin{array}{ll}
\hskip45pt s(h)\otimes \Q&=~h_*({\mathcal L}(N) \cap [N]_{\Q}) -
{\mathcal L}(M) \cap [M]_{\Q}\\[1ex]
&\in  \S_n(M)\otimes \Q~=~
\dot H_n(M;\L_{\bullet})\otimes \Q~=~\sum\limits_{4k \neq n}H_{n-4k}(M;\Q)~.
\hskip65pt\square
\end{array}$$

\noindent{\it Example.} Smale \cite{Smale} proved the generalized
Poincar\'e conjecture: if $N$ is a differentiable $n$-dimensional
manifold with a homotopy equivalence $h:N \to S^n$ and $n \geqslant 5$ then
$h$ is homotopic to a homeomorphism.  Stallings and Newman then proved
the topological version: if $N$ is a topological $n$-dimensional
manifold with a homotopy equivalence $h:N \to S^n$ and $n \geqslant 5$ then
$h$ is homotopic to a homeomorphism.  This is the geometric content of
the computation of the structure set of $S^n$
$$\S^{TOP}(S^n)~=~\S_{n+1}(S^n)~=~0~~(n \geqslant 5)~.\eqno{\square}$$

Here are three consequences of the Theorem in the non-simply-connected case,
subject to the canonical restriction $n \geqslant 5$ : 
\begin{itemize} 
\item[(i)] For any finitely presented group $\pi$ the image of the
assembly map 
$$A~:~H_n(K(\pi,1);\L_{\bullet}) \to L_n(\Z[\pi])$$ 
is the subgroup consisting of the surgery obstructions $\sigma_*(f,b)$
of normal maps $(f,b):N \to M$ of closed $n$-dimensional manifolds with
$\pi_1(M)=\pi$.
\item[(ii)] The Novikov conjecture for a group $\pi$ is that
the higher signatures for any manifold $M$ with $\pi_1(M)=\pi$ 
$$\sigma_x(M)~=~\langle x \cup {\mathcal L}(M),[M] \rangle \in \Q~~
(x \in H^*(K(\pi,1);\Q))$$
are homotopy invariant. The conjecture holds for $\pi$ if and only
if the rational assembly maps 
$$A~:~H_n(K(\pi,1);\L_{\bullet})\otimes \Q~=~H_{n-4*}(K(\pi,1);\Q)
\to L_n(\Z[\pi])\otimes\Q$$
are injective.
\item[(iii)] The topological Borel rigidity conjecture 
for an $n$-dimensional aspherical manifold $M=K(\pi,1)$ is that
every simple homotopy equivalence of manifolds $h:N \to M$ is
homotopic to a homeomeorphism, i.e. $\S^{TOP}(M)=\{*\}$, 
and more generally that
$$\S^{TOP}_{\partial}(M \times D^k,M \times S^{k-1})~=~\{*\}~~(k \geqslant 1)~.$$ 
The conjecture holds for $\pi$ if and only if the
assembly map 
$$A~:~H_{n+k}(K(\pi,1);\L_{\bullet}) \to L_{n+k}(\Z[\pi])$$ 
is injective for $k=0$ and an isomorphism for $k \geqslant 1$.
\end{itemize}

See Chapter 23 of Ranicki \cite{Ranicki1992} and Chapter 8 of Ranicki
\cite{Ranicki1995} for the algebraic Poincar\'e transversality
treatment of the splitting obstruction theory for homotopy equivalences
of manifolds along codimension $q$ submanifolds, involving natural
morphisms $\S_*(X) \to LS_{*-q-1}$ to the $LS$-groups defined
geometrically in Chapter 11 of Wall \cite{Wall}.  The case $q=1$ is
particularly important~: a homotopy invariant functor is a
homology theory if and only if it has excision, and excision is
a codimension 1 transversality property.

\subsection{Homology manifolds} 

An $n$-dimensional Poincar\'e complex $X$ has a {\it 4-periodic total
surgery obstruction} $\overline{s}(X) \in \overline{\S}_n(X)$ such that
$\overline{s}(X)=0$ if (and for $n \geqslant 6$ only if) $X$ is simple
homotopy equivalent to a compact $ANR$ homology manifold (Bryant,
Ferry, Mio and Weinberger \cite{BFMW}).  
The $\S$- and $\overline{\S}$-groups are related by an exact sequence
$$0 \to \S_{n+1}(X) \to \overline{\S}_{n+1}(X) \to
H_n(X;L_0(\Z)) \to \S_n(X) \to \overline{\S}_n(X) \to 0~.$$
The total surgery obstruction $s(X)\in \S_n(X)$ of an $n$-dimensional
homology manifold $X$ is the image of the Quinn \cite{Quinn2}
resolution obstruction
$i(X) \in H_n(X;L_0(\Z))$, such that $i(X)=0$ if (and for $n \geqslant 6$ only if)
there exists a map $M \to X$ from an $n$-dimensional topological manifold $M$
with contractible point inverses. The homology manifold surgery sequence of
$X$ with $n \geqslant 6$ is in one-one correspondence with a portion of the
4-periodic algebraic surgery exact sequence, by a bijection
$$\xymatrix@R+5pt@C+5pt{
\dots \ar[r]& L_{n+1}(\Z[\pi_1(X)]) \ar@{=}[d] \ar[r] &
\S^H(X) \ar[d]^{\displaystyle{\overline{s}}}_-{\cong} \ar[r] & 
[X,L_0(\Z)\times G/TOP]
\ar[d]^{\displaystyle{\overline{t}}}_-{\cong} \ar[r] & 
L_n(\Z[\pi_1(X)]) \ar@{=}[d]\\
\dots \ar[r]& L_{n+1}(\Z[\pi_1(X)]) \ar[r] &\overline{\S}_{n+1}(X) \ar[r] &
H_n(X;\overline{\L}_{\bullet}) \ar[r]^-{\displaystyle{A}} &
L_n(\Z[\pi_1(X)])}$$
with $\S^H(X)$ the structure set of simple homotopy equivalences $h:Y
\to X$ of $n$-dimensional homology manifolds, up to $s$-cobordism.  
\medskip

\noindent{\it Example.} The homology manifold structure set of $S^n$
($n \geqslant 6$) is
$$\S^H(S^n)~=~\overline{\S}_{n+1}(S^n)~=~L_0(\Z)~,$$
detected by the resolution obstruction.\hfill\qed
\medskip

See Chapter 25 of Ranicki \cite{Ranicki1992} and Johnston and Ranicki
\cite{JohnstonRanicki} for more detailed accounts of
the algebraic surgery classification of homology manifolds.
\medskip

The homology manifold surgery exact sequence of \cite{BFMW} 
required the controlled algebraic surgery exact sequence
$$H_{n+1}(B;\overline{\L}_{\bullet}) \to \S_{\epsilon,\delta}(N,f) \to
[N,\partial N;G/TOP,*] \to H_n(B;\overline{\L}_{\bullet})$$
which has now been established by Pedersen, Quinn and Ranicki \cite{PQR}.


\newpage
\addcontentsline{toc}{section}{References}
\providecommand{\bysame}{\leavevmode\hbox to3em{\hrulefill}\thinspace}

\end{document}